\documentclass[aps,preprintnumbers,nofootinbibt,referee]{revtex4}
\usepackage{makeidx}
\usepackage{amsfonts}
\usepackage{amsmath}
\usepackage{amssymb}
\usepackage{eurosym}
\usepackage{graphicx}

\def\be{\begin{equation}}
\def\ee{\end{equation}}
\def\bea{\begin{eqnarray}}
\def\eea{\end{eqnarray}}

\begin{document}

\title{Kosambi-Cartan-Chern (KCC) theory for higher order dynamical systems}
\author{Tiberiu Harko}
\email{t.harko@ucl.ac.uk}
\affiliation{Department of Mathematics, University College London, Gower Street, London
WC1E 6BT, United Kingdom,}
\author{Praiboon Pantaragphong}
\email{kppraibo@kmitl.ac.th}
\affiliation{Mathematics Department, King Mongkut's Institute of Technology, Ladkrabang
BKK 10520, Thailand,}
\author{Sorin V. Sabau}
\email{sorin@tokai.ac.jp}
\affiliation{School of Science, Department of Mathematics, Tokai University, Sapporo
005-8600, Japan}

\begin{abstract}
The Kosambi-Cartan-Chern (KCC) theory represents a powerful mathematical
method for the investigation of the properties of dynamical systems. The KCC
theory introduces a geometric description of the time evolution of a
dynamical system, with the solution curves of the dynamical system described
by methods inspired by the theory of geodesics in a Finsler spaces. The
evolution of a dynamical system is geometrized by introducing a non-linear
connection, which allows the construction of the KCC covariant derivative,
and of the deviation curvature tensor. In the KCC theory the properties of
any dynamical system are described in terms of five geometrical invariants,
with the second one giving the Jacobi stability of the system. Usually, the
KCC theory is formulated by reducing the dynamical evolution equations to a
set of second order differential equations. In the present paper we
introduce and develop the KCC approach for dynamical systems described by
systems of arbitrary $n$-dimensional first order differential equations. We
investigate in detail the properties of the $n$-dimensional autonomous
dynamical systems, as well as the relationship between the linear stability and the
Jacobi stability. As a main result we find that only
even-dimensional dynamical systems can exhibit both Jacobi stability and instability behaviors,
while odd-dimensional dynamical systems are always Jacobi unstable, no
matter their Lyapunov stability. As applications of the developed
formalism we consider the geometrization and the study of the Jacobi stability of the complex dynamical networks, and of the $\Lambda $-Cold Dark Matter cosmological models, respectively.

\end{abstract}

\maketitle
\date{\today }

\section{Introduction}

Continuously time evolving dynamical systems are one of the fundamental
theoretical tools for the mathematical modelling of the evolution of natural
phenomena. They are widely used in every branch of natural sciences, and
they provide extremely useful and powerful descriptions of the evolution of
various natural processes. We define a deterministic continuous dynamical
systems as a set of formal mathematical rules that describe the evolution of
points (elements) included in a set $S$ with respect to an external
continuous time parameter $t\in T$ \cite{Punzi}. Even that generally one can
take the time parameter as being discrete, or a continuous function \cite%
{Punzi}, in the following we restrict our analysis to the case of a
continuous time parameter. More exactly, from a formal mathematical point of
view, we  define a solution of a dynamical system as a map  \cite{Punzi}
\begin{equation}
\phi:T \times S \rightarrow S, (t,x)\mapsto \phi (t,x),
\end{equation}
satisfying the condition $\phi (t , \cdot) \circ \phi (s , \cdot)=\phi (t+s
, \cdot)$, $\forall t ,s\in T$. Usually, in the case of realistic dynamical
systems that may be required to model natural phenomena, one must add some
additional mathematical conditions to the above definition.

In most of the applications in natural sciences the usefulness of the
dynamical systems is largely determined by their predictive capabilities. On
the other hand, the predictive potential of the dynamical systems is
determined essentially by the stability of their solutions. We must also
note that in realistic scientific investigations some uncertainties in the
measured initial conditions in physical, chemical or biological systems do
always appear, and they must be taken into account when formulating the
basic mathematical equations of the considered processes. Therefore, in
order for a dynamical system to be useful in practical applications, it must
provide essential information, and control, not only on the time evolution
the system, but also on the possible deviations of the trajectories of the
dynamical system with respect to a given trajectory, taken as reference. It is important to point out that from the study of a
dynamical system we must obtain two essential information: (i) the local
description of the stability, and (ii) the global description of the
late-time deviations of the trajectories.

The global stability of the solutions of the dynamical systems is described,
by using a rigorous mathematical approach, by the well established and
intensely investigated theory of the Lyapunov stability. In this standard
approach to the problem of stability, the fundamental mathematical
quantities of interest are the Lyapunov exponents. The Lyapunov exponents
measure the exponential deviations of the trajectories of a dynamical system
from a given reference trajectory \cite{1,2}. However, in trying to
determine the Lyapunov exponents one should be aware that usually it is very
difficult to obtain them in an exact analytic form. Therefore, in order to
estimate the Lyapunov exponents one must use numerical methods. Presently
many numerical methods have been developed for their computation, and they
are used extensively in the analysis of the time behavior of the dynamical
systems \cite{3,12}.

The methods of the Lyapunov stability analysis are well established, and
they offer a clear understanding of some stability properties of the
dynamical systems. However, in order to give a deeper insight into behavior
of natural systems, alternative methods for the study of the stability must
also be considered, and developed. Then, once such a new stability method is
obtained, one can compare in detail the predictions and results of the new
proposed method with the results obtained by using the Lyapunov linear
stability analysis of the given dynamical system.

One of the approaches that could offer important new insights for the study
of the stability of the dynamical systems is represented by what one may
call the geometro-dynamical analysis. One of the first examples of such an
alternative study of the stability is the Kosambi-Cartan-Chern (KCC) theory.
The early formulation of the KCC theory was introduced in the important
works of Kosambi \cite{Ko33}, Cartan \cite{Ca33} and Chern \cite{Ch39},
respectively. From a mathematical point of view the KCC theory is inspired
by the geometry of the Finsler spaces, which represents its theoretical
foundation. The KCC theory is based on the fundamental assumption that there
is a one to one correspondence between a second order, autonomous or
non-autonomous, dynamical system, and the geodesic equations in a Finsler
space that can always be associated to the given system (for a recent review
of the KCC theory see \cite{rev}). Essentially, the KCC theory represents a
differential geometric approach to the variational differential equations
describing the deviations of the entire trajectory of a dynamical system
with respect to the nearby ones \cite{An00}. In the geometrical description
introduced by the KCC theory one can associate to each dynamical system two
connections, the first one being a non-linear connection, while the second
one is a Berwald type connection. Five geometrical invariants can be
constructed with the help of the non-linear and Berwald connections,
respectively. The most important of these five geometrical invariants is the
second invariant, called the curvature deviation tensor. Its importance from
the point of view of scientific and mathematical applications relies on the
fact that it gives the so-called Jacobi stability of the given dynamical
system \cite{rev, An00, Sa05,Sa05a}. Different physical, biochemical, or
engineering systems have been extensively studied with the use of the KCC
theory \cite{Sa05, Sa05a, An93, YaNa07, Ha1, Ha2, T0, T1, Ab1, Ab2, Ab3,
Ab4, Ha3, Ha4}.

The study of the two-dimensional autonomous systems, and of their stability
properties, represents an important field of application of the KCC theory.
In \cite{Sa05a, Sa05} two dimensional systems of the form
\begin{equation}  \label{4n}
\frac{du}{dt}=f(u,v),\qquad \frac{dv}{dt}=g(u,v),
\end{equation}
were studied under the assumption that the point $(0,0)$ is a fixed point,
i.e. $f(0,0)=g(0,0)=0$. In the approach introduced in \cite{Sa05a, Sa05},
after relabeling $v$ as $x$, and $g(u,v)$ as $y$, and by assuming that $%
g_u|_{(0,0)}\neq 0$, it follows that one can eliminate the variable $u$.
Since $(u,v)=(0,0)$ is a fixed point, it follows from the Theorem of
Implicit Functions that in the vicinity of $(x,y)=(0,0)$ the equation $%
g(u,x)-y=0$ has a solution $u=u(x,y)$ . Since $\ddot x = \dot g = g_u \, f +
g_v \, y$, an autonomous one-dimensional second order equation is obtained,
which is equivalent to the system (\ref{4n}), and which is given by
\begin{equation}  \label{5n}
\ddot x^1 + g^1(x,y) = 0,
\end{equation}
where
\begin{equation}
g^1(x,y)=-g_u(u(x,y),x) \, f(u(x,y),x) - g_v(u(x,y),x) \, y.
\end{equation}
Hence the Jacobi stability properties of Eq.~(\ref{5n}) can be studied in
detail by using the KCC theory \cite{Sa05a,Sa05}. Consequently, the KCC
stability properties of the first order system (\ref{4n}) can also be
obtained, and the in the in depth comparison between the Jacobi and Lyapunov
stability properties of the two dimensional dynamical systems can be easily
performed.

Usually, the KCC theory is applied to dynamical systems by reducing the
evolution equations to a set of second order differential equations. An
alternative view of the KCC theory, in which the theory is applied to a
first order dynamical system, was introduced in \cite{Ha4}. The properties
of the two dimensional autonomous dynamical systems were investigated in
detail, and the relationship between the linear stability and the Jacobi
stability was also established. As a physical application of the proposed
approach the KCC type geometrization of Hamiltonian systems with one degree
of freedom was considered, and the KCC stability properties of Hamiltonian
systems were discussed.

It is the goal of the present paper to generalize and extend the perspective
introduced in \cite{Ha4} for arbitrary $n$-dimensional first order dynamical
systems. By starting with an $n$-dimensional autonomous first order
dynamical system, we transform it "directly" to an equivalent second order
autonomous system of ordinary differential equations. This transformation is
done by simply taking the time derivative of all equations of the dynamical
system with respect to the time parameter. The obtained second order
equations can be investigated with the use of the KCC theory, which allows a
full geometric description of the system. In particular, the expressions of
the non-linear connection, of the Berwald connection, and of the deviation
curvature tensor are explicitly obtained. Hence an in depth comparison
between the Lyapunov stability and the KCC stability of the critical points
can be performed. In particular, we obtain the important result that only
\textit{even dimensional} first order dynamical systems exhibit Jacobi
stability. We also discuss in detail the case of the three-dimensional
dynamical systems, which are always Jacobi unstable, no matter their
Lyapunov stability.

As a first application of the obtained results we briefly point out the
possibility of investigating the KCC stability of complex dynamic networks.
The theory and applications of complex networks represents a very active
field of research in science and technology. In particular, collective
motions of coupled dynamical networks, such as stabilization,
synchronization, robustness, diffusion etc., have received a great deal of
attention recently \cite{1}, \cite{net1, net2, net3}. We show that a KCC
type geometric description of complex network dynamics is always possible,
and that this description allows the study of the Jacobi stability of the
network.

The theoretical results obtained in the study of the dynamical systems have
also found important applications in the field of cosmology \cite%
{cosm1,cosm2,cosm3,cosm4}. The Lyapunov stability properties of some
cosmological models have been considered in detail \cite{cosm3,cosm4}. The
investigation of the critical points in the phase space for given
cosmological models can provide very important information on the overall
cosmological dynamics. The existence of cosmological attractors can be
related to cosmological solutions indicating the early or late-time behavior
of the Universe. The saddle equilibrium points may give some indications on
the transient cosmological solutions \cite{cosm4}. Hence, the study of the
properties of the equilibrium points of the cosmological models can offer
essential information on the model behavior. As a further step in the study
of the cosmological models from a dynamical systems perspective we consider
the KCC geometrization and the Jacobi stability properties of the $\Lambda $%
Cold Dark Matter ($\Lambda $CDM) models. The KCC approach provides the
possibility of a full geometric description of the $\Lambda $CDM
cosmologies, and the study of their Jacobi stability properties.

The present paper is organized as follows. In Section~\ref{sect1} we briefly
review the Lyapunov stability of systems of ordinary differential equations.
A brief presentation of the basic concepts and results of the KCC theory is
presented in Section~\ref{sect2}. The comparison of the Lypunov and Jacobi
stability for $n$-dimensional dynamical systems is performed in Section~\ref%
{sect3}. The case of the three-dimensional dynamical systems is discussed in
Section~\ref{sect4}. The applications of the KCC theory for the
geometrization and stability study of complex networks are presented in
Section~\ref{sect5}. The possibility of the KCC geometrization of the
cosmological dynamics, and the study of the Jacobi stability of the $\Lambda
$CDM model is considered in Section~\ref{cosm}. We discuss and conclude our
results in Section~\ref{sect6}. In the Appendix we briefly discuss the
relation between the stability of the critical points located at zero, and
the arbitrary critical points.

\section{Lyapunov stability of systems of Ordinary Differential Equations}

\label{sect1}

In the present Section we briefly review following \cite{M} some results on
the Lyapunov stability of the dynamical systems, described by systems of
Ordinary Differential Equations (ODEs). We first recall that, in general,
the stability of a system of ODEs is determined by the roots of its
characteristic polynomial. Indeed, let us consider the system of ODEs
\begin{eqnarray}  \label{1}
\frac{dx^{1}}{dt} &=&f^{1}\left( x^{1},x^{2},...,x^{n}\right) ,  \notag \\
\frac{dx^{2}}{dt} &=&f^{2}\left( x^{1},x^{2},...,x^{n}\right) ,  \notag \\
&&..............., \\
\frac{dx^{n-1}}{dt} &=&f^{n-1}\left( x^{1},x^{2},...,x^{n}\right) ,  \notag
\\
\frac{dx^{n}}{dt} &=&f^{n}\left( x^{1},x^{2},...,x^{n}\right) ,  \notag
\end{eqnarray}
where $f_1, f_2, \dots , f_n$ are smooth functions. Let us linearize the above system about a steady state $\left(
x_{0}^{1},x_{0}^{2},...,x_{0}^{n}\right) $ by associating the linear system
\begin{equation}
\begin{pmatrix}
\frac{dx^{1}}{dt} \\
\frac{dx^{2}}{dt} \\
... \\
\frac{dx^{n}}{dt}%
\end{pmatrix}%
=A%
\begin{pmatrix}
x^{1} \\
x^{2} \\
... \\
x^{n}%
\end{pmatrix}%
,  \label{2}
\end{equation}%
where $A=\left. J\left( f^{1},f^{2},...,f^{n}\right) \right\vert _{\left(
x_{0}^{1},x_{0}^{2},...,x_{0}^{n}\right) }$ is the Jacobian matrix of the
system (\ref{1}), evaluated at the steady state,
\begin{equation*}
A=\left.
\begin{pmatrix}
\frac{\partial f^{1}}{\partial x^{1}} & \frac{\partial f^{1}}{\partial x^{2}}
& ... & \frac{\partial f1}{\partial x^{n}} \\
... & ... & ... & ... \\
\frac{\partial f^{n}}{\partial x^{1}} & \frac{\partial f^{n}}{\partial x^{2}}
& ... & \frac{\partial f^{n}}{\partial x^{n}}%
\end{pmatrix}%
\right\vert _{\left( x_{0}^{1},x_{0}^{2},...,x_{0}^{n}\right) }.
\end{equation*}

The solutions of Eq. (\ref{2}) are given by \cite{M}
\begin{equation}
\begin{pmatrix}
x^{1} \\
x^{2} \\
... \\
x^{n}%
\end{pmatrix}%
=%
\begin{pmatrix}
C^{1}e^{\lambda _{1}t} \\
C^{2}e^{\lambda _{2}t} \\
... \\
C^{n}e^{\lambda _{n}t}%
\end{pmatrix}%
,  \label{3}
\end{equation}%
where $\left( C^{1},C^{2},...,C^{n}\right) $ is a constant vector and $%
\left( \lambda _{1},\lambda _{2},...,\lambda _{n}\right) $ are the
eigenvalues of the matrix $A$, obtained as the roots of the characteristic
polynomial $p(\lambda)$ given by
\begin{equation}
\det \left( A-\lambda I_{n}\right) =0,
\end{equation}%
where $I_{n}$ is the identity matrix .

\textbf{Definition}. \textit{A solution of the dynamical system (\ref{1}) is
called stable if and only if all the roots $\lambda _{1}$, $\lambda _{2}$,
..., $\lambda _{n}$, of the characteristic polynomial $p(\lambda) $ lie in
the left hand side of the complex plane, i. e., $\mathrm{Re}\;\lambda <0$
for all roots $\lambda $.}

If this is true, then $x^{i}(t)=e^{\lambda _{i}t}\rightarrow 0$
exponentially as $t\rightarrow \infty $ for all $i$, and hence $\left(
x^{1},x^{2},...,x^{n}\right) =\left( 0,0,...,0\right) $ is stable to small
(linear) perturbations.

In the $n$-dimensional case, the characteristic polynomial reads
\begin{equation}
\mathbf{p}(\lambda )=\lambda ^{n}+a_{1}\lambda ^{n-1}+...+a_{n}=0,
\label{pol}
\end{equation}%
where the coefficients $a_{i}$, $i=1,2,...,n$ are real numbers. Without any
loose of generality we can always assume $a_{n}\neq 0$, because otherwise we
have $\lambda =0$ and we obtain a characteristic polynomial of order $n-1$
with the zero order coefficient non-vanishing.

The necessary and sufficient conditions for $\mathbf{p}\left( \lambda
\right) $ to have all solutions such that $\mathrm{Re}\lambda <0$ can be
written as \cite{M}
\begin{equation}
a_{n}>0, D_{1}=a_{1}>0, D_{2}=%
\begin{vmatrix}
a_{1} & a_{3} \\
1 & a_{2}%
\end{vmatrix}%
>0, D_{3}=%
\begin{vmatrix}
a_{1} & a_{3} & a_{5} \\
1 & a_{2} & a_{4} \\
0 & a_{1} & a_{3}%
\end{vmatrix}%
>0, ...,
\end{equation}

\begin{equation}
D_{k}=%
\begin{vmatrix}
a_{1} & a_{3} & ... & ... \\
1 & a_{2} & a_{4} & ... \\
0 & a_{1} & a_{3} & ... \\
0 & 1 & a_{2} & ... \\
... & ... & ... & ... \\
0 & 0 & ... & a_{k}%
\end{vmatrix}%
>0,
\end{equation}%
for all $k=1,2,...,n$.

Another method to study the stability problem is to consider the relations
between the non-vanishing solutions of $\mathbf{p}\left( \lambda \right) $.
We have
\begin{equation}
s:=\sum_{i=1}^{n}\lambda _{i}=-a_{1},
\end{equation}
\begin{equation}
\mu _{1}:=\sum_{i,j,i\neq j}^{n}\lambda _{i}\lambda _{j}=a_{2},
\end{equation}
\begin{equation}
...
\end{equation}
\begin{equation}
p:=\lambda _{1}...\lambda _{n}=\left( -1\right) ^{n}a_{n}.
\end{equation}

By taking into account the values of these coefficients we can obtain
important information about the stability of the system (\ref{1}).

\textbf{Remark} (Descartes' Rule of Signs) \cite{M}

Consider the polynomial (\ref{pol}) with $a_{n}>0$. Let $m$ be the number of
sign changes in the sequence of coefficients $\left\{
a_{n},a_{n-1},...,a_{0}\right\} $, ignoring any coefficients which are zero.
Then there are at most $m$ roots of $\mathbf{p}(\lambda )$ which are real
and positive, and further there $m$, $m-2$, $m-4$, ..., real positive roots
\cite{M}. By setting $\omega :=-\lambda $ this property gives important
information about the possible real negative roots, which is essential for
the study of stability.

\section{KCC theory and Jacobi stability of dynamical systems}

\label{sect2}

In the present Section we briefly outline the fundamental concepts, basic
assumptions, and results of the KCC theory. We also introduce the relevant
notations, and present the definitions of the relevant geometric objects
(for a detailed presentation see \cite{rev} and \cite{An00}).

\subsection{The geometry of arbitrary dynamical systems}

In the following we assume that the dynamical variable $x^i$, $i=1,2,...,n$
are defined on a real, smooth $n$-dimensional manifold $\mathcal{M}$. The
tangent bundle of $\mathcal{M}$ is denoted by $T\mathcal{M}$. Usually $%
\mathcal{M}$ is thought as $R^n$, $\mathcal{M}=R^n$, and consequently $T%
\mathcal{M}=TR^n=R^n$. Let's consider a subset $\Omega $ of the Euclidian $%
(2n+1)$ dimensional space $R^{n}\times R^{n}\times R^{1}$. On $\Omega $ we
introduce a $2n+1$ dimensional coordinate system $\left(x^i,y^i,t\right)$, $%
i=1,2,...,n$, where $\left( x^{i}\right) =\left(
x^{1},x^{2},...,x^{n}\right) $, $\left( y^{i}\right) =\left(
y^{1},y^{2},...,y^{n}\right) $ and $t$ is the usual time coordinate. We
define the coordinates $y^i$ as
\begin{equation}
y^{i}=\left( \frac{dx^{1}}{dt},\frac{dx^{2}}{dt},...,\frac{dx^{n}}{dt}%
\right) .
\end{equation}

A basic assumption in the KCC theory is that the time coordinate $t$ is an
absolute invariant. Therefore, within the framework of our approach the only
admissible coordinate transformations are of the form
\begin{equation}
\tilde{t}=t,\tilde{x}^{i}=\tilde{x}^{i}\left( x^{1},x^{2},...,x^{n}\right)
,i\in \left\{1 ,2,...,n\right\} .  \label{ct}
\end{equation}

In many situations of scientific interest the equations of motion of a
dynamical system can be derived from a Lagrangian function $L:T\mathcal{M}%
\rightarrow R$ via the Euler-Lagrange equations,
\begin{equation}
\frac{d}{dt}\frac{\partial L}{\partial y^{i}}-\frac{\partial L}{\partial
x^{i}}=F_{i},i=1,2,...,n,  \label{EL}
\end{equation}%
where in the case of mechanical systems $F_{i}$, $i=1,2,...,n$, represents
the external force. If the Lagrangian $L$ is regular, by means of simple
computations one can show that the Euler-Lagrange equations Eq.~(\ref{EL})
are equivalent to a system of second-order ordinary differential equations
given by
\begin{equation}
\frac{d^{2}x^{i}}{dt^{2}}+2G^{i}\left( x^{j},y^{j},t\right) =0,i\in \left\{
1,2,...,n\right\} ,  \label{EM}
\end{equation}%
where each function $G^{i}\left( x^{j},y^{j},t\right) $ is $C^{\infty }$ in
a neighborhood of some initial conditions $\left( x _{0},y _{0},t_{0}\right)
$ in $\Omega $.

%The system given by
%Eq.~(\ref{EM}) is equivalent to a vector field (semispray) $S$, where
%\begin{equation}
%S=y^{i}\frac{\partial }{\partial x^{i}}-2G^{i}\left( x^{j},y^{j},t\right)
%\frac{\partial }{\partial y^{i}}.
%\end{equation}%

The basic idea of the KCC theory is that if an arbitrary system of
second-order differential equations of the form (\ref{EM}) is given, with no
\textit{a priori} Lagrangean function known, we can still study the behavior
of its trajectories by using differential geometric methods. This study can
be performed by using the close analogy with the trajectories of the
standard Euler-Lagrange system.

\subsection{The non-linear and Berwald connections, the curvature deviation
tensor, and the KCC invariants associated to a dynamical system}

In order to study of the geometrical properties associated to the dynamical
system defined by Eqs.~(\ref{EM}), on the base manifold $\mathcal{M}$ we
introduce first the nonlinear connection $N$, with coefficients $N_{j}^{i}$,
defined as \cite{MHSS}
\begin{equation}  \label{ncon}
N_{j}^{i}=\frac{\partial G^{i}}{\partial y^{j}}.
\end{equation}

From a geometric point of view the nonlinear connection $N_{j}^{i}$ can be
interpreted in terms of a dynamical covariant derivative $\nabla ^N$ as
follows. Given two vector fields $v$, $w$ defined over a manifold $\mathcal{M%
}$, we introduce the covariant derivative $\nabla ^N$ as \cite{Punzi}
\begin{equation}  \label{con}
\nabla _v^Nw=\left[v^j\frac{\partial }{\partial x^j}w^i+N^i_j(x,y)w^j\right]%
\frac{\partial }{\partial x^i}.
\end{equation}

For $N_i^j(x,y)=\Gamma _{il}^j(x)y^l$, from Eq.~(\ref{con}) we immediately
recover the definition of the covariant derivative for the special case of a
standard linear connection, as usually defined in Riemmannian geometry.

On the open subset $\Omega \subseteq R^{n}\times R^{n}\times R^{1}$ on which
the dynamical system is defined, and for the non-singular coordinate
transformations defined by Eqs.~(\ref{ct}), we introduce the KCC-covariant
differential of a vector field $\xi ^{i}(x)$ through the definition \cite%
{An93,An00,Sa05,Sa05a}
\begin{equation}  \label{KCC}
\frac{D\xi ^{i}}{dt}=\frac{d\xi ^{i}}{dt}+N_{j}^{i}\xi ^{j}.
\end{equation}

For $\xi ^{i}=y^{i}$ we obtain
\begin{equation}
\frac{Dy^{i}}{dt}=N_{j}^{i}y^{j}-2G^{i}=-\epsilon ^{i}.
\end{equation}
$\epsilon ^{i}$, representing a contravariant vector field, defined on the
subset $\Omega $ of the Euclidian space, is called the first KCC invariant.
From a physical point of view, and in the framework of classical mechanics,
the first KCC invariant $\epsilon ^{i}$ may be interpreted as an external
force.

As a next step in our study we vary the trajectories $x^{i}(t)$ of the
system (\ref{EM}) into nearby ones according to the rules
\begin{equation}  \label{var}
\tilde{x}^{i}\left( t\right) =x^{i}(t)+\eta \xi ^{i}(t), \tilde{y}^{i}\left(
t\right) =y^{i}(t)+\eta \frac{d\xi ^{i}(t)}{dt},
\end{equation}
where $\left| \eta \right|<<1 $ is a small parameter, and $\xi ^{i}(t)$ are
the components of a contravariant vector field, which is defined along the
trajectory $x^{i}(t)$ of the dynamical system. By substituting Eqs.~(\ref%
{var}) into Eqs.~(\ref{EM}), and by taking the limit $\eta \rightarrow 0$,
we obtain the so-called deviation, or Jacobi, equations as given by \cite%
{An93,An00,Sa05,Sa05a}
\begin{equation}  \label{def}
\frac{d^{2}\xi ^{i}}{dt^{2}}+2N_{j}^{i}\frac{d\xi ^{j}}{dt}+2\frac{\partial
G^{i}}{\partial x^{j}}\xi ^{j}=0.
\end{equation}

By using the KCC-covariant derivative, defined by Eq.~(\ref{KCC}), we can
rewrite Eq.~(\ref{def}) in a covariant form as
\begin{equation}
\frac{D^{2}\xi ^{i}}{dt^{2}}=P_{j}^{i}\xi ^{j},  \label{JE}
\end{equation}
where we have denoted
\begin{equation}  \label{Pij}
P_{j}^{i}=-2\frac{\partial G^{i}}{\partial x^{j}}-2G^{l}G_{jl}^{i}+ y^{l}%
\frac{\partial N_{j}^{i}}{\partial x^{l}}+N_{l}^{i}N_{j}^{l}+\frac{\partial
N_{j}^{i}}{\partial t}.
\end{equation}
In Eq.~(\ref{Pij}) we have also introduced the tensor $G_{jl}^{i}$, defined
as \cite{rev, An00, An93,MHSS,Sa05,Sa05a}
\begin{equation}
G_{jl}^{i}\equiv \frac{\partial N_{j}^{i}}{\partial y^{l}},
\end{equation}
and which is called the Berwald connection.

The tensor $P_{j}^{i}$ is the second KCC-invariant, and it may be called
alternatively as the deviation curvature tensor. We shall call Eq.~(\ref{JE}%
) the Jacobi equation. In either Riemann or Finsler geometry, when the
system of equations (\ref{EM}) describes the geodesic motion, then Eq.~(\ref%
{JE}) is the Jacobi field equation corresponding to the given geometry.

A scalar quantity which can be constructed from $P^i_j$ is the trace $P$ of
the curvature deviation tensor, which can be obtained from the relation
\begin{equation}
P=P_{i}^{i}=-2\frac{\partial G^{i}}{\partial x^{i}}-2G^{l}G_{il}^{i}+ y^{l}%
\frac{\partial N_{i}^{i}}{\partial x^{l}}+N_{l}^{i}N_{i}^{l}+\frac{\partial
N_{i}^{i}}{\partial t}.
\end{equation}

In the KCC theory one can also introduce the third, fourth and fifth
invariants of the second order system of equations (\ref{EM}). These
invariants are defined according to \cite{An00}
\begin{equation}  \label{31}
P_{jk}^{i}\equiv \frac{1}{3}\left( \frac{\partial P_{j}^{i}}{\partial y^{k}}-%
\frac{\partial P_{k}^{i}}{\partial y^{j}}\right) ,P_{jkl}^{i}\equiv \frac{%
\partial P_{jk}^{i}}{\partial y^{l}},D_{jkl}^{i}\equiv \frac{\partial
G_{jk}^{i}}{\partial y^{l}}.
\end{equation}

From a geometrical point of view the third KCC invariant $P_{jk}^{i}$ can be
interpreted as a torsion tensor. The fourth and fifth KCC invariants $%
P_{jkl}^{i}$ and $D_{jkl}^{i}$ represent the Riemann-Christoffel curvature
tensor, and the Douglas tensor, respectively \cite{rev, An00}. It is
important to point out that in a Berwald space these tensors always exist.
In the KCC theory these five invariants are the basic mathematical
quantities describing the geometrical properties and interpretation of an
arbitrary system of second-order differential equations.

\subsection{Jacobi stability of dynamical systems}

In many scientific applications, including the study of the stability of
physical, chemical, biological or engineering systems, the study of the
behavior of the trajectories of the dynamical system given by Eqs.~(\ref{EM}%
) in the vicinity of a point $x^{i}\left( t_{0}\right) $ is of extreme
importance. In the following for simplicity we take the origin of the time
coordinate at $t_{0}=0$. We interpret the trajectories $x^{i}=x^{i}(t)$ of
the system of equations (\ref{EM}) as curves in the Euclidean space $\left(
R^{n},\left\langle .,.\right\rangle \right) $, where $\left\langle
.,.\right\rangle $ represents the canonical inner product of $R^{n}$. As for
the deviation vector $\xi $ we assume that it obeys the set of initial
conditions $\xi \left( 0\right) =O$ and $\dot{\xi}\left( 0\right) =W\neq O$,
respectively, where $O\in R^{n}$ is the null vector \cite{rev, An00,
Sa05,Sa05a}.

To describe the focusing/dispersing tendency of the trajectories around $%
t_{0}=0$ we introduce the following mathematical description: if the
condition $\left| \left| \xi \left( t\right) \right| \right| <t^{2}$, $%
t\approx 0^{+}$ is satisfied, then the trajectories are bunching together.
On the other hand, if the deviation vector satisfies the condition $\left|
\left| \xi \left( t\right) \right| \right| >t^{2}$, $t\approx 0^{+}$, it
follows that the trajectories of the system (\ref{EM}) have a dispersing
behavior \cite{rev, An00, Sa05,Sa05a}. The focusing/dispersing behavior of
the trajectories of a dynamical system can be also described from a
geometric point of view by using the properties of the deviation curvature
tensor as follows. For $t\approx 0^{+}$ the trajectories of the system of
equations Eqs.~(\ref{EM}) are bunching/focusing together if and only if the
real parts of the eigenvalues (characteristic values) of the deviation
tensor $P_{j}^{i}\left( 0\right) $ are strictly negative. On the other hand,
if and only if the real part of the eigenvalues (characteristic values) of $%
P_{j}^{i}\left( 0\right) $ are strictly positive, the trajectories are
dispersing \cite{rev, An00, Sa05,Sa05a}.

Based on the above qualitative discussion, we introduce the rigorous
definition of the concept of Jacobi stability for a general continuously
time evolving dynamical system as follows \cite{rev, An00,Sa05,Sa05a}:

\textbf{Definition:} \textit{Let's assume that with respect to the norm $%
\left| \left| .\right| \right| $ induced by a positive definite inner
product, the system of differential equations Eqs.~(\ref{EM}) satisfies the
initial conditions
\begin{equation*}
\left| \left| x^{i}\left( t_{0}\right) -\tilde{x}^{i}\left( t_{0}\right)
\right| \right| =0, \left| \left| \dot{x}^{i}\left( t_{0}\right) -\tilde{x}%
^{i}\left( t_{0}\right) \right| \right| \neq 0.
\end{equation*}
}

\textit{Then the trajectories of the dynamical system given by Eqs.~(\ref{EM}%
) are called Jacobi stable if and only if the real parts of the
characteristic values of the deviation tensor $P_{j}^{i}$ are strictly
negative everywhere. If the real parts of the characteristic values of the
deviation tensor $P_{j}^{i}$ are strictly positive everywhere, the
trajectories of the dynamical system are called Jacobi unstable.}

%Graphically, the focussing/dispersing behavior of the trajectories near the
%origin is depicted in Fig.~\ref{pict1}.
%\begin{widetext}
%\begin{centering}
%\begin{figure*}[htp]
%\centering
%\includegraphics[height=6cm,width=12cm]{fig1_pict.eps} %\centering
%\setlength{\unitlength}{1cm}
%\begin{picture}(10,5)\label{fig:Figure1}
%\thicklines \qbezier(0,1)(2,4)(5,4) \qbezier(0,1)(3,3)(5,1)
%\put(3,1.98){\vector(0,1){1.65}} \put(3,3.63){\vector(0,-1){1.65}}
%\put(4,1){$x^i(t)$} \put(4,4.2){$\tilde{x}^i(t)$}
%\put(3.2,2.8){$\eta(t)$} \put(1,0){$||\xi (t)||^2<t^2,\ t\approx
%0^+$} \put(-0.2,0.5){$x^i(0)$} \qbezier(6,2.5)(8.5,2.5)(10.5,4)
%\qbezier(6,2.5)(8.5,2.5)(10.5,1)
%\put(9.3,1.75){\vector(0,1){1.54}} \put(9.5,1){$x^i(t)$}
%\put(9.5,4.2){$\tilde{x}^i(t)$} \put(7,0){$||\xi (t)||^2>t^2,\
%t\approx 0^+$} \put(5.7,2.0){$x^i(0)$}
%\end{picture}
%%%%%%%%%%%%%%%%%%%%%%%%%%%%%%%%%%%%%%%%%%%%%%%%%%%%%%%%
%\caption{Behavior of the trajectories near zero.}
%\label{pict1}
%\end{figure*}
%\end{centering}
%\end{widetext}

\section{Relation of Jacobi and Lyapunov stability for higher dimensional
dynamical systems}

\label{sect3}

In the following we will always assume, without loosing the generality, that
the fixed point of the system of equations Eqs.~(\ref{1}) is $\left(
0,0,...,0\right) $ (see the Appendix for a discussion of the case of
arbitrary fixed points). Let us consider the system of ODEs given by Eqs.~(%
\ref{1}). We point first out that, geometrically, we can describe the
solutions of Eqs.~(\ref{1}) as a flow $\varphi _{t}:D\subset
R^{n}\rightarrow R^{n}$, or, more generally, $\varphi _{t}:D\subset \mathcal{%
M}\rightarrow \mathcal{M}$, where $\mathcal{M}$ is a smooth $n$-dimensional
manifold. The canonical lift of $\varphi _{t}$ to the tangent space $T%
\mathcal{M}$ is given by
\begin{equation}
\hat{\varphi}_{t}:T\mathcal{M}\rightarrow T\mathcal{M}, \hat{\varphi}%
_{t}\left( u\right) =\left( \varphi _{t}\left( u\right) ,\dot{\varphi}%
_{t}\left( u\right) \right) ,
\end{equation}
where
\begin{equation}
\dot{\varphi}_{t}\left( u\right) =\frac{\partial \varphi \left( t,u\right) }{%
\partial t}.
\end{equation}

\subsection{Geometrical description of higher dimensional dynamical systems}

In terms of dynamical systems we simply take the derivative of (\ref{1})
with respect to the time parameter $t$, thus obtaining
\begin{equation*}
\frac{d^{2}x^{1}}{dt^{2}}=f_{1}^{1}\left( x^{1},...,x^{n}\right)
y^{1}+f_{2}^{1}\left( x^{1},...,x^{n}\right) y^{2}+...+f_{n}^{1}\left(
x^{1},...,x^{n}\right) y^{n},
\end{equation*}
\begin{equation*}
\frac{d^{2}x^{2}}{dt^{2}}=f_{1}^{2}\left( x^{1},...,x^{n}\right)
y^{1}+f_{2}^{2}\left( x^{1},...,x^{n}\right) y^{2}+...+f_{n}^{2}\left(
x^{1},...,x^{n}\right) y^{n},
\end{equation*}
\begin{equation*}
.....................,
\end{equation*}
\begin{equation*}
\frac{d^{2}x^{n}}{dt^{2}}=f_{1}^{n}\left( x^{1},...,x^{n}\right)
y^{1}+f_{2}^{n}\left( x^{1},...,x^{n}\right) y^{2}+...+f_{n}^{n}\left(
x^{1},...,x^{n}\right) y^{n},
\end{equation*}%
where
\begin{equation}
f_{j}^{i}=\frac{\partial f^{i}\left( x^{1},...,x^{n}\right) }{\partial x^{j}}%
,i,j=1,...,n,y^{i}=\frac{dx^{i}}{dt},i=1,...,n.
\end{equation}

In other words, on $T\mathcal{M}$ we obtain
\begin{equation*}
\frac{dy^{1}}{dt}=\sum _{i=1}^{n}f_{i}^{1}\left( x^{1},...,x^{n}\right)
y^{i},
\end{equation*}
\begin{equation*}
.....................,
\end{equation*}
\begin{equation*}
\frac{dy^{n}}{dt}=\sum _{i=1}^{n}f_{i}^{n}\left( x^{1},...,x^{n}\right)
y^{i},
\end{equation*}
where $\left( x^{1},x^{2},...,y^{1},...,y^{n}\right) $ are local coordinates
on $T\mathcal{M}$. Hence one can see that the above system is actually a
linear dynamical system on the fiber $T_{\left( x^{1},x^{2},...,x^{n}\right)
}\mathcal{M}$. The system can be written as
\begin{equation}  \label{29n}
\frac{d^{2}x^{i}}{dt^{2}}-\sum _{k=1}^{n}f_{k}^{i}\left(
x^{1},..,x^{n}\right) y^{k}=0,i=1,...,n.
\end{equation}

By comparing Eq.~(\ref{29n}) to Eq.~(\ref{EM}) we find
\begin{equation}
G^{i}=-\frac{1}{2}\sum_{k=1}^{n}f_{k}^{i}\left( x^{1},..,x^{n}\right)
y^{k},i=1,...,n,  \label{Gi}
\end{equation}%
or, equivalently,
\begin{equation}
\begin{pmatrix}
G^{1} \\
G^{2} \\
. \\
. \\
. \\
G^{n}%
\end{pmatrix}%
=-\frac{1}{2}%
\begin{pmatrix}
\sum_{k=1}^{n}f_{k}^{1}\left( x^{1},..,x^{n}\right) y^{k} \\
\sum_{k=1}^{n}f_{k}^{2}\left( x^{1},..,x^{n}\right) y^{k} \\
. \\
. \\
. \\
\sum_{k=1}^{n}f_{k}^{n}\left( x^{1},..,x^{n}\right) y^{k}%
\end{pmatrix}%
=-\frac{1}{2}J\left( f^{1},f^{2},...,f^{n}\right)
\begin{pmatrix}
y^{1} \\
y^{2} \\
. \\
. \\
. \\
y^{n}%
\end{pmatrix}%
=-\frac{1}{2}J\cdot y.
\end{equation}

Using Eq. (\ref{ncon}) we obtain
\begin{equation}
N_{j}^{i}=-\frac{1}{2}\frac{\partial f^{i}\left( x^{1},..,x^{n}\right) }{%
\partial x^{j}}=-\frac{1}{2}f_{j}^{i}\left( x^{1},..,x^{n}\right)
,i,j=1,...,n,  \label{Nij}
\end{equation}%
or
\begin{equation}
\left( N_{j}^{i}\right) _{i,j=1,...,n}=%
\begin{pmatrix}
N_{1}^{1} & N_{2}^{1} & ... & N_{n}^{1} \\
. & . & . & . \\
N_{1}^{n} & N_{2}^{n} & ... & N_{n}^{n}%
\end{pmatrix}%
=-\frac{1}{2}J\left( f^{1},f^{2},...,f^{n}\right) ,
\end{equation}%
and hence
\begin{equation*}
G_{jl}^{i}=\frac{\partial N_{j}^{i}}{\partial y^{l}}=0.
\end{equation*}

Then, for the components of the deviation curvature tensor $\left(
P_{j}^{i}\right) $, given by Eq.~(\ref{Pij}), we find
\begin{equation}
P_{j}^{i}=-2\frac{\partial G^{i}}{\partial x^{j}}+y^{k}\frac{\partial
N_{j}^{i}}{\partial x^{k}}+N_{l}^{i}N_{j}^{l}.
\end{equation}

Then, with the use of Eqs. (\ref{Gi}) and (\ref{Nij}) we obtain
\begin{equation}
P_{j}^{i}=\frac{1}{2}\sum_{k=1}^{n}f_{jk}^{i}\left( x^{1},..,x^{n}\right)
y^{k}+\frac{1}{4}\sum_{l=1}^{n}f_{l}^{i}\left( x^{1},..,x^{n}\right)
f_{j}^{l}\left( x^{1},..,x^{n}\right) .
\end{equation}

Therefore we have obtained the following

\textbf{Proposition 1}. \textit{The curvature deviation tensor associated to
an $n$-dimensional dynamical system is given by
\begin{equation}
P=\frac{1}{2}\sum_{k,l=1}^{n}%
\begin{pmatrix}
H_{f^{k}}\cdot y & H_{f^{l}}\cdot y%
\end{pmatrix}%
^{t}+\frac{1}{4}J^{2}\left( f^{1},f^{2},...,f^{n}\right) ,
\end{equation}%
where $H_{f^{k}}=%
\begin{pmatrix}
f_{11}^{k} & f_{12}^{k} & ... & f_{1n}^{k} \\
... & ... & ... & ... \\
f_{n1}^{k} & f_{n2}^{k} & ... & f_{nn}^{k}%
\end{pmatrix}%
$ is the Hessian of $f^{k}\left( x^{1},x^{2},...,x^{n}\right) $.}

\subsection{Lyapunov stability versus Jacobi stability}

Evaluating $P_{j}^{i}$ at the stability state $\left( 0,...,0\right) $ we
obtain

\begin{equation}
\left.\left(P_{j}^{i}\right)\right |_{\left( 0,...,0\right) }=-2\left.\left(%
\frac{\partial G^{i}}{\partial x^{j}}\right)\right |_{\left( 0,...,0\right)
}+\left.\left( N_{l}^{i}N_{j}^{l}\right) \right |_{\left( 0,...,0\right) }=%
\frac{1}{4}J^{2}|_{\left( 0,...,0\right) }=\frac{1}{4}A^{2},
\end{equation}%
where $A$ is the Jacobian matrix, with coefficients evaluated at the fixed
point. We briefly discus the case when the stability state is not at zero in
the Appendix.

Let us now recall some general results from linear algebra.

\textbf{Lemma 1}: Let $A$ be an $\left( n\times n\right) $ matrix, and
denote by $\lambda _{1}$, $\lambda _{2}$,..., $\lambda _{n}$ its
eigenvalues. Then

(i) The eigenvalues of $kA$ are $k\lambda _{1}$,...,$k\lambda _{n}$, for any
scalar $k\neq 0$.

(ii) The eigenvalues of $A^{k}:=\underset{k\;\mathrm{times}}{\underbrace{%
A\times ...\times A}}$ are $\left( \lambda _{1}\right) ^{k},...,\left(
\lambda _{n}\right) ^{k}$.

\textbf{Proof}. Let us recall that if $A$ is a complex matrix with complex
eigenvalues $\lambda _1,..., \lambda _n$, counted with their complex
multiplicities, then always there exists an invertible matrix $Q$ such that
the matrix $T:=QAQ^{-1}$ is upper triangular, and the entries on the
diagonal are exactly $\lambda _1, ...,\lambda _n$. This means that $A$ and $%
T $ have the same characteristic polynomial $\mathbf{p}_A=\mathbf{p}%
_T=\left(\lambda -\lambda _1\right)...\left(\lambda -\lambda _n\right)=0$
(see for instance the book \cite{M}). Using this observation the proof of
the \textbf{Lemma} is simple. Indeed, remark that $T=QAQ^{-1}$ implies $%
kT=Q\left(kA\right)Q^{-1}$, where $k$ is a scalar. Then, remarking that
diagonal entries in $kT$ are $k\lambda _1,...,k\lambda _n$, we get $\mathbf{p%
}_A=\mathbf{p}_T=\left(\lambda -k\lambda _1\right)\left(\lambda -k\lambda
_2\right)...\left(\lambda -k\lambda _n\right)$, and therefore (i) is proved.

Similarly, by remarking that $T=QAQ^{-1}$ implies $T^k=Q\cdot A^k\cdot
Q^{-1} $, and that the diagonal entries of $T^k$ are $\left(\lambda
_1\right)^k,...,\left(\lambda _k\right)^k$, then (ii) follows by the same
arguments as above.

From \textbf{Lemma 1} we obtain an important result that relates the two
types of stability.

\textbf{Lemma 2}: If $\lambda _{1}$, $\lambda _{2}$,..., $\lambda _{n}$ are
eigenvalues of the matrix $A=J\left( f_{1},f_{2},...,f_{n}\right) |_{\left(
0,...,0\right) }$ then $\mu _{1}=\left( \lambda _{1}\right) ^{2}$,..., $\mu
_{n}=\left( \lambda _{n}\right) ^{2}$ are eigenvalues of the deviation
curvature matrix $P|_{\left( 0,...,0\right) }=P_{j}^{i}|_{\left(
0,...,0\right) }$, $i,j=1,...,n$.

This \textbf{Lemma} leads to a very natural relation between classical
stability theory and KCC theory.

\textbf{Theorem 1}. \textit{Let us consider the ODEs system (\ref{1}), and
denote by $A:=J\left( f_{1},f_{2},...,f_{n}\right) |_{\left( 0,...,0\right)
} $ the Jacobian matrix evaluated at the fixed point $\left( 0,...,0\right) $%
. }

\textit{1. If $A$ has only real eigenvalues, that is, the steady state is a
(stable or unstable) node, then $P|_{\left( 0,...,0\right) }$ has only real
positive eigenvalues, that is, the steady state is always Jacobi unstable. }

\textit{2. If $A$ has conjugate complex eigenvalues $\mu _{j}:=\alpha
_{j}+i\beta _{j}$, $\bar{\mu}:=\alpha _{j}-i\beta _{j}$ for some $j$ so that
$2j\leq n$, then }

\textit{2.1. If $2j=n$, then the steady state Jacobi stable if and only if $%
\alpha _{j}^{2}-\beta _{j}^{2}<0$, for all $j=1,...,n/2$; }

\textit{2.2. If $2j<n$, then $p\left( \lambda \right) $ has real and complex
eigenvalues and therefore the steady state is Jacobi unstable. }

We can now summarize our main results as follows.

\textbf{Theorem 2}. \textit{The steady state $\left( 0,...,0\right) $ of the
dynamical system (\ref{1}) is Jacobi stable if and only if (\ref{1})
satisfies the following }

\textit{(i) It is even dimensional, i.e. $n=2j$. }

\textit{(ii) All eigenvalues of the characteristic polynomial $p(\lambda )=0$
are complex conjugate, i.e. $\lambda _{j}:=\alpha _{j}+i\beta _{j}$, $\bar{%
\lambda}_{j}:=\alpha _{j}-i\beta _{j}$, $j=1,...,n/2$. }

\textit{(iii) For any $j=1,...,n/2$ we have
\begin{equation*}
\alpha _{j}^{2}-\beta _{j}^{2}<0.
\end{equation*}
}

\section{The case of the three variable dependent dynamical system}

\label{sect4}

In the following we consider the dynamical system
\begin{equation*}
\frac{dx^{1}}{dt}=f\left( x^{1},x^{2},x^{3}\right) ,
\end{equation*}
\begin{equation}
\frac{dx^{2}}{dt}=g\left( x^{1},x^{2},x^{3}\right) ,  \label{53}
\end{equation}
\begin{equation*}
\frac{dx^{3}}{dt}=h\left( x^{1},x^{2},x^{3}\right) ,
\end{equation*}%
and consider again that the steady state $\left(
x_{0}^{1},x_{0}^{2},x_{0}^{3}\right) =\left( 0,0,0\right) $. The
characteristic polynomial associated to the above system reads
\begin{equation}
\mathbf{p}\left( \lambda \right) =\lambda ^{3}-s\cdot \lambda ^{2}+\mu \cdot
\lambda -p,
\end{equation}%
where we have defined
\begin{equation}
s:=\lambda _{1}+\lambda _{2}+\lambda _{3},
\end{equation}
\begin{equation}
\mu :=\lambda _{1}\lambda _{2}+\lambda _{2}\lambda _{3}+\lambda _{3}\lambda
_{1},
\end{equation}
\begin{equation}
p=\lambda _{1}\lambda _{2}\lambda _{3}.
\end{equation}

A simple computation of the determinant
\begin{equation}
\left\vert A-\lambda I_{3}\right\vert =0,
\end{equation}%
where
\begin{equation}
A=%
\begin{pmatrix}
a_{11} & a_{12} & a_{13} \\
a_{21} & a_{22} & a_{23} \\
a_{31} & a_{32} & a_{33}%
\end{pmatrix}%
=\left.
\begin{pmatrix}
\frac{\partial f\left( x^{1},x^{2},x^{3}\right) }{\partial x^{1}} & \frac{%
\partial f\left( x^{1},x^{2},x^{3}\right) }{\partial x^{2}} & \frac{\partial
f\left( x^{1},x^{2},x^{3}\right) }{\partial x^{3}} \\
\frac{\partial g\left( x^{1},x^{2},x^{3}\right) }{\partial x^{1}} & \frac{%
\partial g\left( x^{1},x^{2},x^{3}\right) }{\partial x^{2}} & \frac{\partial
g\left( x^{1},x^{2},x^{3}\right) }{\partial x^{3}} \\
\frac{\partial h\left( x^{1},x^{2},x^{3}\right) }{\partial x^{1}} & \frac{%
\partial h\left( x^{1},x^{2},x^{3}\right) }{\partial x^{2}} & \frac{\partial
h\left( x^{1},x^{2},x^{3}\right) }{\partial x^{3}}%
\end{pmatrix}%
\right\vert _{\left( 0,0,0\right) },
\end{equation}%
and $I_{3}$ is the $\left( 3,3\right) $ identity matrix, implies the
familiar form
\begin{equation}
\mathbf{p}\left( \lambda \right) =\lambda ^{3}-\mathrm{tr}\;A\cdot \lambda
^{2}+\frac{1}{2}\left[ \left( \mathrm{tr}\;A\right) ^{2}-\mathrm{tr}\left(
A^{2}\right) \right] \cdot \lambda -\det A=0,
\end{equation}%
where $\mathrm{tr}\;A$ and $\det A$ are the trace and determinant of the
matrix $A$, and $\mathrm{tr}\left( A^{2}\right) $ is the trace of the matrix
$A^{2}:=A\cdot A$. The classification of the steady states of a
three-dimensional dynamical system was given by Poincar\'{e}, and can be
summarized as follows:

1. $\lambda _{1}$, $\lambda _{2}$, $\lambda _{3}$ real roots of same sign of
$\mathbf{p}\left( \lambda \right) $: $p$ is called a node.

2. $\lambda _{1}$, $\lambda _{2}$, $\lambda _{3}$ real roots of $\mathbf{p}%
\left( \lambda \right) $, but not of same sign: $p$ is called a saddle.

3. $\lambda _{1},$ $\lambda _{2}\in \mathbb{C}$, $\lambda _{3}\in \mathbb{R}$%
, such that $\mathrm{Re}\;\left( \lambda _{1,2}\right) $ and $\lambda _{3}$
are of the same sign: $p$ is called a focus.

4. $\lambda _{1},$ $\lambda _{2}\in \mathbb{C}$, $\lambda _{3}\in \mathbb{R}$%
, such that $\mathrm{Re}\;\left( \lambda _{1,2}\right) $ and $\lambda _{3}$
are not of the same sign: $p$ is called a saddle-focus.

5. $\lambda _{1},$ $\lambda _{2}\in \mathbb{C}$, $\lambda _{3}\in \mathbb{R}$%
, such that $\mathrm{Re}\;\left( \lambda _{1,2}\right) =0$: $p$ is called a
center.

This classification can be combined with the Lyapunov stability in the case $%
\mathrm{Re}\;\left( \lambda _{i}\right) <0$ for $i=1,2,3$. Taking now
\textbf{Theorem 2} into account it follows

\textbf{Theorem 3}. \textit{The fixed point $p$ of the three-dimensional
dynamical system }(\ref{53}) \textit{is always classified as Jacobi
unstable, regardless its Lyapunov stability}.

\textbf{Remark}. Our \textbf{Theorem 3} above might be understood that in
the case $n=3$, and more generally $n=2j$, $j\in
%TCIMACRO{\U{2115} }%
%BeginExpansion
\mathbb{N}
%EndExpansion
$, KCC theory gives no useful information on the steady states. This is
completely true. Without having any Jacobi stable trajectory, if $p$ is a
focus or a saddle-focus, then the case when $\lambda _{1},$ $\lambda _{2}\in
%TCIMACRO{\U{2102} }%
%BeginExpansion
\mathbb{C}
%EndExpansion
$, $\lambda _{3}\in
%TCIMACRO{\U{211d} }%
%BeginExpansion
\mathbb{R}
%EndExpansion
$, $\alpha ^{2}-\beta ^{2}<0$, where $\lambda _{1,2}=\alpha \pm i\beta $,
gives information about a Jacobi-type saddle-focus, distinct from the usual
Jacobi unstable behavior. We will not pursue here further investigations of
these behaviors.

\section{Jacobi stability analysis of dynamical networks}

\label{sect5}

The study of the dynamics of complex
networks is has many applications in both science and engineering. Different processes like such as
stabilization, synchronization, robustness, diffusion etc. can be described by using networks of coupled systems. From mathematical point of view a complex network of coupled systems is described by the system of ODEs \cite{1}
\begin{equation}
\frac{dx^{i}}{dt}=F^{i}\left( x^{1},x^{2},...,x^{N}\right) -\sigma \sum_{r}{%
L_{r}^{i}H^{r}\left( x^{1},x^{2},...,x^{N}\right) }%
=f^{i}(x^{1},x^{2},...,x^{N}),i=1,2,...,N,  \label{1n}
\end{equation}%
where $x^{i}$ with $i\in {1,2,...,N}$ are dynamical variables, $F$ and $H$
are evolution and coupling functions, respectively, $\sigma $ is a constant,
and $L_{r}^{i}$ is the Laplacian matrix, defined by $L_{i}^{i}=k^{i}$ (the
connectivity degree of node $i$), $L_{r}^{i}=-1$ if nodes $i$ and $j$ are
connected, and $L_{r}^{i}=0$ otherwise.

A Lyapunov linear stability analysis of the dynamical system describing networks of coupled systems can be performed by (i) expanding
around a critical state $x_1 =x_2 = ...=x_N = x^s_{cr}$, with $x^s_{cr}$ solution
of $\dot{x}^ s_{cr} = F\left(x^s_{cr}\right)=0$, (ii) diagonalizing $L$ to find its $N$
eigenvalues $0 =\lambda _1 <\lambda _2\leq ...\leq \lambda _N$, and, finally,  (iii)
writing the governing equations for the normal modes $y_i$ of the perturbations
\begin{equation}
\dot{y}_i=\left[F^{\prime }\left(x^s_{cr}\right)-\sigma \lambda _iH^{\prime
}\left(x^s_{cr}\right)\right]y_i.
\end{equation}
The above equations have all the same form, but different effective couplings $\alpha
=\sigma \lambda _i$.

As a first step in the Jacobi stability analysis of the network dynamics we
need to extend the dynamical process given by Eq.~(\ref{1n}) to
differentiable manifolds. By taking the derivative of Eqs.~(\ref{1n}) with
respect to the time we obtain
\begin{equation}  \label{28}
\frac{d^2x^i}{dt^2}=\frac{\partial f^i}{\partial x^l}y^l=\left(\frac{%
\partial F^i}{\partial x^l}y^l-\sigma L^i_r\frac{\partial H^r}{\partial x^n}%
y^n\right)=\left(\frac{\partial F^i}{\partial x^n}-\sigma L^i_r\frac{%
\partial H^r}{\partial x^n}\right)y^n.
\end{equation}

Eq.~(\ref{28}) has the same mathematical form as Eq.~(\ref{EM}), with
\begin{equation}
G^{i}\left( x^{a},y^{b}\right) =-\frac{1}{2}\left( \frac{\partial F^{i}}{%
\partial x^{n}}-\sigma L_{r}^{i}\frac{\partial H^{r}}{\partial x^{n}}\right)
y^{n}.
\end{equation}%
Therefore we can apply now the KCC theory for the description of the
dynamics of the network system. The non-linear connection associated to the
system Eq.~(\ref{28}) is given by
\begin{equation}
N_{j}^{i}=\frac{\partial G^{i}}{\partial y^{j}}=-\frac{1}{2}\left( \frac{%
\partial F^{i}}{\partial x^{n}}-\sigma L_{r}^{i}\frac{\partial H^{r}}{%
\partial x^{n}}\right) \delta _{j}^{n}=-\frac{1}{2}\left( \frac{\partial
F^{i}}{\partial x^{j}}-\sigma L_{r}^{i}\frac{\partial H^{r}}{\partial x^{j}}%
\right) .
\end{equation}

All the components of the Berwald connection $G_{jl}^{i}=\partial
N_{j}^{i}/\partial y^{l}$ are identically equal to zero. Then we easily
obtain
\begin{equation}
-2\frac{\partial G^{i}}{\partial x^{j}}=\left( \frac{\partial ^{2}F^{i}}{%
\partial x^{j}\partial x^{n}}-\sigma L_{r}^{i}\frac{\partial ^{2}H^{r}}{%
\partial x^{j}\partial x^{n}}\right) y^{n},y^{n}\frac{\partial N_{j}^{i}}{%
\partial x^{n}}=-\frac{1}{2}\left( \frac{\partial ^{2}F^{i}}{\partial
x^{j}\partial x^{n}}-\sigma L_{r}^{i}\frac{\partial ^{2}H^{r}}{\partial
x^{j}\partial x^{n}}\right) y^{n},
\end{equation}%
giving for the deviation curvature tensor the expression
\begin{equation}
P_{j}^{i}=\frac{1}{2}\left( \frac{\partial ^{2}F^{i}}{\partial x^{j}\partial
x^{n}}-\sigma L_{r}^{i}\frac{\partial ^{2}H^{r}}{\partial x^{j}\partial x^{n}%
}\right) y^{n}+\frac{1}{4}\left( \frac{\partial F^{i}}{\partial x^{l}}%
-\sigma L_{r}^{i}\frac{\partial H^{r}}{\partial x^{l}}\right) \left( \frac{%
\partial F^{l}}{\partial x^{j}}-\sigma L_{r}^{l}\frac{\partial H^{r}}{%
\partial x^{j}}\right) .
\end{equation}

The trace of the deviation curvature tensor of the complex network is given
by
\begin{equation}
P=P_{i}^{i}=\frac{1}{2}\left( \frac{\partial ^{2}F^{i}}{\partial
x^{i}\partial x^{n}}-\sigma L_{r}^{i}\frac{\partial ^{2}H^{r}}{\partial
x^{i}\partial x^{n}}\right) y^{n}+\frac{1}{4}\left( \frac{\partial F^{i}}{%
\partial x^{l}}-\sigma L_{r}^{i}\frac{\partial H^{r}}{\partial x^{l}}\right)
\left( \frac{\partial F^{l}}{\partial x^{i}}-\sigma L_{r}^{l}\frac{\partial
H^{r}}{\partial x^{i}}\right) .
\end{equation}

Therefore a dynamical network is Jacobi stable if and only if the real parts
of the eigenvalues of the deviation tensor $P_{j}^{i}$ are strictly negative
everywhere, and Jacobi unstable otherwise. By assuming that the critical
point of the network is given by $(0,...,0)$, it follows that the deviation
curvature tensor evaluated at the critical point has the form
\begin{equation}
\left.P_j^i\right|_{(0,...,0)}=\left. \frac{1}{4}\left( \frac{\partial F^{i}%
}{\partial x^{l}}-\sigma L_{r}^{i}\frac{\partial H^{r}}{\partial x^{l}}%
\right) \left( \frac{\partial F^{l}}{\partial x^{j}}-\sigma L_{r}^{l}\frac{%
\partial H^{r}}{\partial x^{j}}\right)\right |_{(0,...,0)}.
\end{equation}

Therefore all our previous results, summarized in \textbf{Theorem 2}, can be
applied for the case of complex networks. In particular it follows that only
even-dimensional networks are Jacobi stable. The geodesic deviation equation
can also be obtained, and its study may give some insights into the
development of chaos in complex networks.

\section{KCC geometrization and the Jacobi stability of the $\Lambda $CDM
cosmological models}

\label{cosm}

In the following we restrict our analysis to isotropic and homogeneous
cosmological models, described by the Friedmann-Robertson-Walker line
element, given by
\begin{equation}
ds^2=-dt^2+a^2(t)\left(dx^2+dy^2+dz^2\right),
\end{equation}
where $a(t)$ is the scale factor of the Universe. The cosmological expansion
for a Universe filled with pressureless dust and radiation is described by
the Friedmann equations, which take the form
\begin{equation}  \label{Fr1}
3H^2=\rho _m+\rho _r+\Lambda, \dot{H}=-\frac{1}{2}\left(\rho _m+\frac{4}{3}%
\rho _r\right),
\end{equation}
\begin{equation}  \label{Fr2}
\dot{\rho }_m+3H\rho _m=0,\dot{\rho}_r+4H\rho _r=0
\end{equation}
where a dot denotes the derivative with respect to the time $t$. In Eqs.~(%
\ref{Fr1}) and (\ref{Fr2}) $H=\dot{a}/a$ is the Hubble function, $\rho _m$
is the matter energy density, $\rho _r$ is the energy density of the
radiation, and $\Lambda $ is the cosmological constant.

\subsection{Friedmann equations as an autonomous dynamical system}

In order to formulate the cosmological evolution in terms of a dynamical
system, we introduce first the matter, radiation and cosmological constant
density parameters $\left( \Omega _{m},\Omega _{r},\Omega _{\Lambda }\right)
$, defined as
\begin{equation}
\Omega _{m}=\frac{\rho _{m}}{3H^{2}},\Omega _{r}=\frac{\rho _{r}}{3H^{2}}%
,\Omega _{\Lambda }=\frac{\Lambda }{3H^{2}}.
\end{equation}%
The density parameters satisfy the relation
\begin{equation}
\Omega _{m}+\Omega _{r}+\Omega _{\Lambda }=1.
\end{equation}%
As the basic variables $(x,y)$ in the phase space we choose $x\equiv \Omega
_{r}$, and $y\equiv \Omega _{\Lambda }$, respectively. Then $\Omega
_{m}=1-x-y$, and the range of the variables is since $0\leq x\leq 1$, $0\leq
y\leq 1$, and $0\leq \Omega _{m}\leq 1$. The physically relevant phase space
for the cosmological dynamics is defined as $\Phi =\left\{ (x,y):x+y\leq
1,0\leq x\leq 1,0\leq y\leq 1\right\} $. By taking the time derivatives of $%
x $ and $y$ with respect to the time, after introducing the new time
variable $\tau =\ln a(t)$, it follows that the Friedmann equations can be
formulated in terms of an autonomous dynamical system given by \cite{cosm3}
\begin{equation}
\frac{dx}{d\tau }=-x(1-x+3y),  \label{ca1}
\end{equation}%
\begin{equation}
\frac{dy}{d\tau }=(3+x-3y)y.  \label{ca2}
\end{equation}%
The critical points of the system (\ref{ca1}) and (\ref{ca2}) in $\Psi $ can
be obtained by solving the algebraic equations $x(1-x+3y)=0$ and $%
(3+x-3y)y=0 $, respectively, and are given by
\begin{equation}
P_{dS}=\{x=0,y= 1\},P_r=\{x= 1,y= 0\},P_m=\{x =0,y= 0\}.
\end{equation}

The Lyapunov stability properties of the system can be obtained from the
study of the Jacobian matrix
\begin{equation}
J=%
\begin{pmatrix}
-1+2x-3y & -3x \\
y & 3+x-6y%
\end{pmatrix}%
.
\end{equation}

The critical point $P_{dS} = (0, 1)$, $\Omega _{\Lambda} = 1$, with $a(t)
\propto \exp\left( \sqrt{\Lambda/3}t\right)$, is associated to an
accelerated de Sitter type expansion. It is a future attractor \cite{cosm4}.
The critical point $P_r = (1, 0)$, with $\Omega _r = 1$ and $a(t)\propto
\sqrt{t}$, corresponds to the radiation-dominated era, and is a source
point, or a past attractor. Finally, the critical point $P_m (0, 0)$, with $%
\Omega _m = 1$ and $a(t) \propto t^{2/3}$, corresponds to the
matter-dominated phase of the cosmological expansion, and to a decelerating
expansion.  $P_m$ is a saddle critical point \cite{cosm4}.

\subsection{The KCC geometrization and the Jacobi stability of the Friedmann
equations}

In order to apply the KCC theory to the dynamical system given by equations (%
\ref{ca1}) and (\ref{ca2}), we relabel the variables as $x\equiv x^1$ and $%
y\equiv x^2$. We also denote $y^1=dx^1/d\tau$ and $y^2=dx^2/d\tau$,
respectively. Hence we obtain
\begin{equation}  \label{ca3}
\frac{dx^1}{d\tau}=-x^1\left(1-x^1+3x^2\right)=f\left(x^1,x^2\right),
\end{equation}
\begin{equation}  \label{ca4}
\frac{dx^2}{d\tau}=\left(3+x^1-3x^2\right)x^2=g\left(x^1,x^2\right).
\end{equation}

Then, by taking the derivative with respect to $\tau $ of Eqs.~(\ref{ca3})
and (\ref{ca4}), we obtain the following lift on the tangent bundle of the
cosmological dynamical system,
\begin{equation}
\frac{d^2x^1}{d\tau ^2}=\left(-1+2x^1-3x^2\right)y^1-3x^1y^2,
\end{equation}
\begin{equation}
\frac{d^2x^2}{d\tau ^2}=x^2y^1+\left(3+x^1-6x^2\right)y^2.
\end{equation}

By comparison with Eqs.~(\ref{EM}) we obtain immediately
\begin{equation}
(G)=G^{i}\left( x^{1},x^{2},y^{1},y^{2}\right) =-\frac{1}{2}%
\begin{pmatrix}
\left( -1+2x^{1}-3x^{2}\right) y^{1}-3x^{1}y^{2} \\
x^{2}y^{1}+\left( 3+x^{1}-6x^{2}\right) y^{2}%
\end{pmatrix}%
.
\end{equation}

The components of the non-linear connection $\left( N\right) =N_{i}^{j}$
associated to the cosmological dynamical system are obtained as
\begin{equation}
\left( N\right) =N_{i}^{j}=-\frac{1}{2}%
\begin{pmatrix}
-1+2x^{1}-3x^{2} & -3x^{1} \\
x^{2} & 3+x^{1}-6x^{2}%
\end{pmatrix}%
.
\end{equation}

For the components of the deviation tensor $P=P_{i}^{j}$ we obtain
\begin{equation}
P=\frac{1}{2}%
\begin{pmatrix}
H_{f}\cdot y & H_{g}\cdot y%
\end{pmatrix}%
^{t}+\frac{1}{4}J^{2}\left( f,g\right) ,
\end{equation}
where $H_{f}=%
\begin{pmatrix}
f_{11} & f_{12} \\
f_{21} & f_{22}%
\end{pmatrix}%
$ is the Hessian of $f$, and $H_g$ is the Hessian of $g$. By evaluating $P$ at the critical points gives
\be
\left.P\right|_{\left(y^1=0,y^2=0,x^1=x^1_{cr},x^2=x^2_{cr}\right)}=\frac{1}{4}A^2,
\ee
where $A=\left.J\right|_{\left(y^1=0,y^2=0,x^1=x^1_{cr},x^2=x^2_{cr}\right)}$. Explicitly, the curvature deviation tensor for the Friedmann cosmological dynamical system can be obtained as
\be
\left.P\right|_{\left(y^1=0,y^2=0\right)}=\left(
\begin{array}{cc}
 \frac{1}{4} (-2 x^1+3 x^2+1)^2+\frac{3 x^1 x^2}{2} & \frac{3}{2} x^1 (x^1-6 x^2+3)+\frac{3}{4} x^1 (-2 x^1+3 x^2+1)
   \\
 (x^1-6 x^2+3) x^2+\frac{1}{2} (-2 x^1+3 x^2+1) x^2 & (x^1-6 x^2+3)^2+\frac{3 x^1 x^2}{2}
\end{array}
\right).
\ee

We compute the numerical values of the curvature deviation tensor first at the critical point $P_{dS}=\left(x^1=0,x^2=1\right)$. In this critical point the curvature deviation tensor takes the form
\be
\left.P\right|_{\left(y^1=0,y^2=0,x^1=0,x^2=1\right)}=
\left(
\begin{array}{cc}
 4 & 0 \\
 -1 & 9
\end{array}
\right)
\ee
and has the eigenvalues $\lambda _1=4$ and $\lambda _2=9$. Hence the critical point $P_{dS}$ of the standard $\Lambda $CDM cosmological model is Jacobi unstable. For the critical point $P_r= \left(x^1=1,x^2=0\right)$, we obtain
\be
\left.P\right|_{\left(y^1=0,y^2=0,x^1=1,x^2=0\right)}=
\left(
\begin{array}{cc}
 \frac{1}{4} & \frac{21}{4} \\
 0 & 16
\end{array}
\right),
\ee
with the corresponding eigenvalues $\lambda _1=16$ and $\lambda _2=1/4$. Hence the critical point $P_r$ is also Jacobi unstable. Finally, for the critical point $P_m=(0,0)$ we obtain
\be
\left.P\right|_{\left(y^1=0,y^2=0,x^1=0,x^2=0\right)}=
\left(
\begin{array}{cc}
 \frac{1}{4} & 0 \\
 0 & 9
\end{array}
\right),
\ee
with the eigenvalues $\lambda _1=9$, $\lambda _2=1/4$. Hence all the critical points of the dynamical system associated to the Friedmann equations are Jacobi unstable.

\section{Discussions and final remarks}

\label{sect6}

In the present paper we have extended and generalized an alternative
approach to the standard KCC theory for first order autonomous dynamical
systems, introduced in \cite{Ha4}. This approach is based on a different,
"direct" transformation of the first order $n$-dimensional system to second
order differential equations. We have presented and discussed this "direct"
approach in detail for arbitrary dimensional dynamical systems. For this
class of arbitrary dimensional dynamical systems we have developed and
discussed in detail the two basic stability analysis methods -- the
(Lyapunov) linear stability analysis and the Jacobi stability analysis. From
the point of view of the KCC theory, such a "direct" approach allows for a
generalization of the geometric framework for first order systems.
Consequently, the parameter space is increased, as well as the predictive
power, of the method. While in the two dimensional case there is a good
correlation between the Lyapunov stability of the critical points, and the
Jacobi stability of the same critical points \cite{Sa05}, a very different
picture emerges for dynamical systems with $n>2$. The fundamental result of
our analysis indicates that only even-dimensional dynamical systems are
Jacobi stable, while odd-dimensional dynamical systems are always Jacobi
unstable, no matter the Lyapunov stability of their critical points.
describing the robustness of the corresponding trajectory to a small
perturbation \cite{Sa05}. On the other hand, it is important to point out
that the Jacobi stability analysis is a very convenient way to describe the
resistance of limit cycles to small arbitrary perturbation of trajectories.

As a first application of our approach we have considered the possibility of a
geometric description of the dynamics of complex networks. The dynamical
behavior of such systems is described by $n$-dimensional dynamical systems,
consisting of systems of first order, highly non-linear, differential
equations. Such systems can be easily interpreted in the geometric framework
of the KCC theory, and the associated geometric quantities (non-linear
connection, Berwald connection, and deviation curvature tensor) can be
obtained easily. The results obtained in this study allow a full geometric
approach to the stability of complex systems, as well as the investigation
of the Jacobi stability of the critical points. Hopefully, future studies
may shed some light on the relation between Jacobi stability and chaotic
behavior in complex networks.

We have also considered the Jacobi stability analysis of the $\Lambda $CDM cosmological models. These models can be naturally represented in terms of an autonomous second degree dynamical system. The KCC theory allows a straightforward geometrization of this system, which can be described in purely geometric  terms with the help of a non-linear connection, the associated covariant derivative, and the curvature deviation tensor. The eigenvalues of the deviation curvature tensor as estimated at the critical points give the Jacobi stability properties of the cosmological model. We have considered three such models, corresponding to the de Sitter expansion, and to the matter and radiation dominated epochs, respectively. Our results show that the critical points corresponding to this cosmological models are all Jacobi unstable.

In summary, in the present paper we have introduced and studied in detail
some geometrical methods necessary for an in depth analysis and description
of the KCC and Lyapunov stability properties of arbitrary dimensional
dynamical systems. Such methods may help scientists to better understand the
time evolution and the dynamical properties of natural phenomena.

\section*{Acknowledgments}

We would like to thank to the anonymous referee for comments and suggestions
that helped us to significantly improve our manuscript.

\appendix
\section{The case of the arbitrary fixed points}

For the sake of simplicity we have assumed in our exposition that $%
O=(0,...,0)$ is the fixed point of ours system of ODEs. This assumption is
not essential, and all the presented results remain true for arbitrary fixed
points. Again for the sake of simplicity we will explain this fact only for
the two-dimensional case. The higher dimensional case can be treated in an
identical way.

Let us consider the autonomous system of ODEs

\begin{equation}
\frac{du}{dt}=f(u,v),\frac{dv}{dt}=g\left( u,v\right) ,  \label{ap1}
\end{equation}%
and consider the variable change

\begin{equation}
\bar{u}:=u-u_{0},\bar{v}=v-v_{0},  \label{ap2}
\end{equation}%
where $u_{0}$, $v_{0}$ are some real constants. For these new variables the
system of ODEs Eqs. (\ref{ap1}) become

\begin{equation}
\frac{d\bar{u}}{dt}=f\left( \bar{u}+u_{0},\bar{v}+v_{0}\right) =:\bar{f}%
\left( \bar{u},\bar{v}\right) ,\frac{d\bar{v}}{dt}=g\left( \bar{u}+u_{0},%
\bar{v}+v_{0}\right) =:\bar{g}\left( \bar{u},\bar{v}\right) .  \label{ap3}
\end{equation}

It is easy to see that ODEs (\ref{ap3}) and (\ref{ap1}) are topologically
and geometrically equivalent, hence they have the same stability behavior at
the fixed points.

Indeed, geometrically we recall that giving (\ref{ap1}) is equivalent in
defining a vector field

\begin{equation}
X|_{(u,v)}:=f(u,v)\frac{\partial }{\partial u}+g(u,v)\frac{\partial }{%
\partial v}.
\end{equation}

A solution of Eqs. (\ref{ap1}) is given by the flow $\varphi ^{t}$ of $X$ in
$\mathbb{R}^{2}$ (or, more generally, in the two-dimensional manifold $%
\mathcal{M}$. Similarly, (\ref{ap3}) gives a vector field $\bar{X}|_{\left(
\bar{u},\bar{v}\right) }$ obtained from $X$ by parallel translation given by
(\ref{ap2}). Since $du/dt=d\bar{u}/dt$, $dv/dt=d\bar{v}/dt$, it is clear
that the direction of $X $ does not change.

The geometrical underlying reason for this result is clear. The group of
isometries of the Euclidian plane $\mathbb{R}^{2}$ is made of translations
and rotations, and hence it is clear that the translation (\ref{ap2})
actually maps any solution curve of (\ref{ap1}) to a solution curve of (\ref%
{ap3}) by keeping the orientation.

From topological point of view, we recall that two flows $\varphi ^{t}$ and $%
\psi ^{t}$ on some manifold $\mathcal{M}$ are topologically equivalent if
there exists a homomorphism $h:\mathcal{M}\longrightarrow \mathcal{M}$ that
takes the trajectories of $\varphi ^{t}$ to the trajectories of $\psi ^{t}$,
while preserving the orientation. In our case, due to the linearity of (\ref%
{ap2}), ODEs (\ref{ap1}) and (\ref{ap3}) are obviously topologically
equivalent. Moreover, they are also geometrically equivalent due to the
discussion above.

\end{document}